\documentclass[final]{siamltex}
\newtheorem{Theorem}{Theorem}[section]
\newtheorem{Lemma}{Lemma}[section]
\newtheorem{Remark}{Remark}[section]
\newcommand{\beq}{\begin{eqnarray}}
\newcommand{\eeq}{\end{eqnarray}}
\newcommand{\beqno}{\begin{eqnarray*}}
\newcommand{\eeqno}{\end{eqnarray*}}
\newcommand{\be}{\begin{equation}}
\newcommand{\ee}{\end{equation}}

\newcommand{\D}{\displaystyle}
\newcommand\qed{\hfill$\Box$}
\usepackage{amssymb,amsmath}

 \allowdisplaybreaks

\title{Global Spherically Symmetric Classical Solution to Compressible Navier-Stokes Equations
with Large Initial Data and Vacuum}

\author{Shijin Ding\thanks{School of Mathematical Sciences, South China Normal University,  Guangzhou  510631,  P.R.  China
({\tt Ding: dingsj@scnu.edu.cn}).} \and Huanyao
Wen$^*$\thanks{Department of Mathematics, Central China Normal
University, Wuhan 430079, P.R. China ({\tt Wen:
huanyaowen@hotmail.com})} \and Lei Yao\thanks{Department of
Mathematics, Northwest University, Xi'an 710127, P.R. China ({\tt
Yao: yaolei1056@hotmail.com}).} \and Changjiang
Zhu\thanks{Department of Mathematics, Central China Normal
University, Wuhan 430079, P.R. China ({\tt Zhu:
cjzhu@mail.ccnu.edu.cn}) (Corresponding author).}}

\begin{document}

\maketitle

\begin{abstract}
 In the paper, we obtain a result on the existence and uniqueness of global
 spherically symmetric classical solutions to the compressible
isentropic Navier-Stokes equations with vacuum in a bounded domain
or exterior domain $\Omega$ of $\mathbb{R}^n$($n\ge2$). Here, the
initial data could be large. Besides, the regularities of the
solutions are better than
 those obtained in \cite{9,10,12}. The analysis is based on some new mathematical techniques and some
new useful energy estimates. This is an extensive work of
\cite{9,10,12}, where the global radially symmetric strong
solutions, the local classical solutions in 3D and the global
classical solutions in 1D were obtained, respectively. This paper
can be viewed as the first result on the existence of global
classical solutions with large initial data and vacuum in higher
dimension.
\end{abstract}



\pagestyle{myheadings} \thispagestyle{plain} \markboth{SHIJIN
DING, HUANYAO WEN, LEI YAO, AND CHANGJIANG ZHU}{COMPRESSIBLE
NAVIER-STOKES EQUATIONS}

\section{Introduction}

In this paper, we consider the initial-boundary value problem of
compressible isentropic Navier-Stokes equations in a bounded
domain or exterior domain $\Omega$ of  $\mathbb{R}^n(n\ge2)$: \beq
\label{1.1}
\begin{cases}
\rho_t+\nabla\cdot(\rho {\bf u})=0,\ \rho\ge0,\\
(\rho {\bf u})_t+\nabla\cdot(\rho {\bf u}\otimes{\bf u})+\nabla
P(\rho)=\mu\triangle{\bf u}+(\mu+\lambda)\nabla (\nabla\cdot{\bf
u})+\rho {\bf f},
\end{cases} \eeq for $({\bf x},t)\in \Omega\times (0,+\infty)$,
where
$$
\Omega=\{{\bf x}\big|a<|{\bf x}|<b\}, 0<a<b\le\infty,\ {\bf
f}({\bf x},t)=f(|{\bf x}|,t)\frac{{\bf x}}{|{\bf x}|},
$$
$\rho$ and $ {\bf u}=(u_1, u_2, \cdots, u_n)$ denote the density
and the velocity respectively; $P(\rho)=K\rho^\gamma$, for some
constants $\gamma>1$ and $K>0$, is the pressure function; ${\bf
f}$ is the external force; the viscosity coefficients $\mu$ and
$\lambda$ satisfy the natural physical restrictions: $\mu>0$ and
$2\mu+n\lambda\ge0$.

We consider the initial condition: \be\label{1.2}(\rho, \ {\bf
u})\big|_{t=0}=(\rho_0,\ {\bf u}_0) \ {\rm{in}}\
\overline{\Omega}, \ee and the boundary condition: \be\label{1.3}
 {\bf u}({\bf x},t)\rightarrow0, \ \textrm{as}\ \ |{\bf x}|\rightarrow a\ \textrm{or}\ b, \ \textrm{for}\ t\geq0,
\ee where
$$
\rho_0({\bf x})=\rho_0(|{\bf x}|),\ {\bf u}_0({\bf x})=u_0(|{\bf
x}|)\frac{{\bf x}}{|{\bf x}|}.
$$
We are looking for a spherically symmetric classical solution
$(\rho,{\bf u})$:
$$
\rho({\bf x},t)=\rho(r,t),\ {\bf u}({\bf x},t)=u(r,t)\frac{{\bf
x}}{r},
$$
where $r=|{\bf x}|$, and $(\rho, u)(r,t)$ satisfies
 \beq
\label{1.4}
\begin{cases}
\displaystyle\rho_t+(\rho u)_r+m\frac{\rho u}{r}=0,\ \rho\ge0,\\
\displaystyle(\rho u)_t+(\rho u^2)_r+m\frac{\rho u^2}{r}+
P_r=\nu(u_r+m\frac{u}{r})_r+\rho f,
\end{cases} \eeq
for $(r,t)\in (a,b)\times(0,\infty)$, with the initial condition:
\be\label{1.5}(\rho(r,t), \ u(r,t))\big|_{t=0}=(\rho_0(r),\
u_0(r)) \ \ {\rm{in}}\ I, \ee and the boundary condition:
\be\label{1.6}
  u(r,t)\rightarrow0, \ \textrm{as}\ r\rightarrow a\ \textrm{or}\ b, \ \textrm{for}\ t\geq0,
\ee where $m=n-1$,  $\nu=2\mu+\lambda\geq \frac{2(n-1)}{n}\mu>0$
and $I=[a, b]$.

Let's review some previous work in this direction. When the
viscosity coefficient $\mu$ is constant, the local classical
solution of non-isentropic Navier-Stokes equations in H\"older
spaces was obtained by Tani in \cite{1} with $\rho_0$ being
positive and essentially bounded. Using delicate energy methods in
Sobolev spaces, Matsumura and Nishida showed in \cite{2,3} that
the existence of the global classical solution   provided that the
initial data was small in some sense and away from vacuum. There
are also some results about the existence of global strong
solution to the Navier-Stokes equations with constant viscosity
coefficient when $\rho_0>0$, refer for instance to \cite{4,5} for
the isentropic flow. Jiang in \cite{jiang1}  proved the global
existence of spherically symmetric smooth solutions in H\"older
spaces to the equations of a viscous polytropic ideal gas in the
domain exterior to a ball in $\mathbb{R}^n$ ($n=2$ or $3$) when
$\rho_0>0$. For general initial data,
 Kawohl in  \cite{6} got the global classical solution
 with $\rho_0>0$ and the viscosity coefficient $\mu=\mu(\rho)$ satisfying
$$
0<\mu_0\le\mu(\rho)\le\mu_1, \ \ {\rm for} \ \ \rho\ge0,
$$
where $\mu_0$ and $\mu_1$ are constants. Indeed, such a condition
includes the case $\mu(\rho)\equiv$const.

In the presence of vacuum, Lions in \cite{lions} used the weak
convergence method  to  show the existence of global weak solution
to the Navier-Stokes equations for isentropic flow with general
initial data and $\gamma\ge\frac{9}{5}$ in  three dimensional
space. Later, the restriction on $\gamma$ was relaxed by Feireisl,
et al \cite{Feireisl} to $\gamma>\frac{3}{2}$. Unfortunately, this
assumption excludes for example the interesting case $\gamma=1.4$
(air, et al). Jiang and Zhang relaxed the condition to $\gamma>1$
in
 \cite{jiang} when they considered the global spherically
symmetric weak solution. It worths mentioning a result due to Hoff
in  \cite{hoff}, who showed the existence of weak solutions for
the case
 $\gamma=1$ with $\rho_0$ being positive and essentially
bounded.

There were few results about strong solution when the initial
density may vanish until  Salvi and Stra$\breve{s}$kraba in
\cite{salvi}, where $\Omega$ is a bounded domain, $P = P(\cdot)\in
C^2[0,\infty)$, $\rho_0\in H^2$, ${\bf u}_0\in H^1_0\bigcap H^2$,
and satisfied the compatibility condition: \be\label{1.7} L{\bf
u}_0({\bf x})-\nabla P(\rho_0)({\bf x})=\sqrt{\rho}_0{\bf g},\
\textrm{for}\ {\bf g}\in L^2,\ee where
$L:=\mu\Delta+(\mu+\lambda)\nabla \mathrm{div}$ is the Lam\'{e}
operator.

 Afterwards, Cho, Choe and Kim in
\cite{7,8,9} established some local and global existence results
about strong solution in   bounded or unbounded domain with
 initial data different  from \cite{salvi} still satisfying
(\ref{1.7}). Particularly, Choe and Kim in \cite{9} showed  the
radially symmetric strong solution existed globally in time for
$\gamma\ge2$ in annular domain. As it is pointed out in \cite{9}
that the results have been proved only for annular domain and
cannot be extended to a ball $\Omega=B_R =\{{\bf x}\in
\mathbb{R}^n: |{\bf x}|<R<\infty\}$, because of a counter-example
of Weigant \cite{wei}. Precisely, for $1<\gamma<1+\frac{1}{n-1}$,
Weigant constructed a radially symmetric strong solution $(\rho,
{\bf u})$ in $B_R\times [0, 1)$ such that $\|\rho(\cdot,
t)\|_{L^\infty(B_R)}\rightarrow \infty$ as $t\rightarrow 1^-$. A
recent paper \cite{11} written by Fan, Jiang and Ni improved the
result in \cite{9} to the case $\gamma\ge1$.

The local classical solution was obtained by Cho and Kim in
\cite{10} when the initial density may vanish and satisfying the
following compatibility conditions:
 \be\label{1.8} L{\bf u}_0({\bf x})-\nabla P(\rho_0)({\bf
x})=\rho_0[{\bf g_1}({\bf x})-{\bf f}({\bf x},0)],\ee for ${\bf
x}\in\overline{\Omega}$, ${\bf g_1}\in D_0^1$ and
$\sqrt{\rho_0}{\bf g_1}\in L^2$. Recently, we  used some new
estimates to get a unique globally classical solution $\rho\in
C^1([0,\infty);H^3)$ and $ u\in H_{loc}^1([0,\infty);H^3)$ to one
dimensional compressible Navier-Stokes equations in a bounded
domain when the initial density may vanish, cf. \cite{12}. It
worths mentioning that Xin in \cite{24} showed that the smooth
solution $(\rho,u)\in C^1([0,\infty);H^3(\mathbb{R}^1))$ to the
Cauchy problem must blow up when the initial density is of
nontrivial compact support.

Since the system (\ref{1.4}) have the one dimensional feature, the
results in \cite{12} are possible to be obtained here. Moreover,
we get higher regularities of the solutions to system (\ref{1.4}).
This causes some new challenges compared with \cite{12}, which
will be handled by some new estimates.

This paper can be viewed to be the first result on global
classical solutions with large initial data and vacuum in higher
dimension.

\bigbreak

 \noindent{\bf
Notations:}

(1) $Q_T=I\times[0,T]$,
$\widetilde{Q}_T=\overline{\Omega}\times[0,T]$ for $T>0$.

(2) For $p\ge 1$, $L^p=L^p(\Omega)$ denotes the $L^p$ space with
the norm $\|\cdot\|_{L^p}$. For $k\ge 1$ and $p\ge 1$,
$W^{k,p}=W^{k,p}(\Omega)$ denotes the Sobolev space, whose norm is
denoted as $\|\cdot\|_{W^{k,p}}$; $H^k=W^{k,2}(\Omega)$.

 (3) For an integer $k\ge 0$ and $0<\alpha<1$, let $C^{k+\alpha}(\overline{\Omega})$ denote
the Schauder space of function on $\overline{\Omega}$, whose $k$th
order derivative is H\"older continuous with exponent $\alpha$,
with the norm $\|\cdot\|_{C^{k+\alpha}}$.

(4) For an integer $k\ge 0$, denote $$
H_r^k=H_r^k(I)\triangleq\left\{u\Big|\sum\limits_{i=0}^k\int_Ir^m|\partial_r^iu|^2<\infty\right\},
$$
with the norm
$$
\|\cdot\|_{H_r^k}=\left(\sum\limits_{i=0}^k\int_Ir^m|\partial_r^iu|^2\right)^\frac{1}{2},
$$
$L_r^2=H_r^0.$

(5) $D^{k,p}=\left\{v\in
L^1_{loc}(\Omega)\big|\|\nabla^kv\|_{L^p}<\infty\right\}$,
$D^k=D^{k,2}$.

(6) $D_0^1=D_0^{1, 2}$  is the closure of $C_0^{\infty}(\Omega)$
in $D^{1, 2}$.

\vspace{4mm}

Our main results are stated as follows.
\begin{Theorem}\label{th:1.1}
Assume that $\rho_0\geq0$ satisfies $\rho_0\in L^1\cap
H^2,\rho_0^\gamma\in H^2$, ${\bf u}_0\in D^3\cap D_0^1$ and ${\bf
f}\in C([0,\infty);H^1)$, ${\bf f}_t\in
L^2_{loc}([0,\infty);L^2)$, and  the initial data $\rho_0, \ {\bf
u}_0$  satisfy the compatible condition (\ref{1.8}) with ${\bf
g}_1({\bf x})=g_1(r)\frac{{\bf x}}{r}$.
 Then for any $T>0$,  there exists a unique global
classical solution $(\rho, {\bf u})$ to (\ref{1.1})-(\ref{1.3})
satisfying \beqno (\rho,\rho^\gamma)\in C([0,T];H^2),\ \rho\geq0,\
 {\bf u} \in C([0,T];D^3\cap D^1_0),\ {\bf u}_t\in
L^\infty([0,T];D_0^1)\cap L^2([0,T];D^2). \eeqno
\end{Theorem}
\begin{Remark}\label{re:1.1}
{\rm{(i) Note that if $\Omega$ is bounded and locally
Lipschitzian, then $D^{k,p}=W^{k,p}$. See \cite{Galdi} for the
proof.

(ii) By Sobolev embedding theorems, we have
$$
H^k(I)\hookrightarrow C^{k-\frac{1}2}(I),\ \ {\rm for} \ \
k=1,2,3,
$$
this together with the regularities of $(\rho, u)$ give
$$(\rho,\rho^\gamma)\in C([0,T];C^{1+\frac{1}{2}}(I)),\ u\in
C([0,T];C^{2+\frac{1}{2}}(I)).
$$
Since $(\rho({\bf x},t),{\bf u}({\bf
x},t))=(\rho(r,t),u(r,t)\frac{{\bf x}}{r})$, we get
$$(\rho,\rho^\gamma)\in C([0,T];C^{1+\frac{1}{2}}(\overline{\Omega})),\ {\bf u}\in
C([0,T];C^{2+\frac{1}{2}}(\overline{\Omega})),
$$
 which means $(\rho,{\bf u})$ is the classical solution to
(\ref{1.1})-(\ref{1.3}).}}
\end{Remark}

\begin{Theorem}\label{th:1.2}
Consider the same assumptions as in Theorem \ref{th:1.1}, and in
addition assume that $\rho_0\in H^5,\rho_0^\gamma\in H^5$,
$\nabla(\sqrt{\rho_0})\in L^\infty$, $\sqrt{\rho_0}\nabla^2{\bf
g_1}\in L^2$, $\rho_0\nabla^3{\bf g_1}\in L^2$, ${\bf u}_0\in
D^5$, ${\bf f}\in C([0,\infty);H^3)\cap
L^2_{loc}([0,\infty);H^4)$, ${\bf f}_t\in C([0,\infty);H^1)\cap
L^2_{loc}([0,\infty);H^2)$ and ${\bf f}_{tt}\in
L^2_{loc}([0,\infty);L^2)$. Then the regularities of the solution
obtained in Theorem \ref{th:1.1} can be improved as follows:
\beqno &(\rho,\rho^\gamma)\in C([0,T];H^5),\ {\bf u}\in
L^\infty([0,T];D_0^1\cap D^5)\cap L^2([0,T];D^6),&\\& {\bf u}_t\in
L^\infty([0,T];D_0^1)\cap L^2([0,T];D^3),\
(\sqrt{\rho}\nabla^2{\bf u}_t, \rho\nabla^3 {\bf u}_t)\in
L^\infty([0,T];L^2).& \eeqno
\end{Theorem}

The rest of the paper is organized as follows. In Section 2, we
prove Theorem \ref{th:1.1}. In Section 3, we prove Theorem
\ref{th:1.2} by giving some estimates similar to \cite{12} and
some new estimates, such as Lemma \ref{le:3.6}, Lemma \ref{le:3.7}
and Lemma \ref{le:3.8}.

The constants $\nu$ and $K$ play no role in the analysis, so we
assume $\nu=K=1$  without loss of generality.
\setcounter{section}{1} \setcounter{equation}{0}
\section{ \ Proof of Theorem \ref{th:1.1}}

In this section, we get a unique global classical solution to
(\ref{1.4})-(\ref{1.6}) with initial density $\rho_0\ge\delta>0$
and $b<\infty$ by some {\it a priori} estimates globally in time
based on the local solution. Moreover, the estimates are
independent of $b$ and $\delta$. Next, we construct a sequence of
approximate solutions to (\ref{1.4})-(\ref{1.6}) under the
assumption $\rho_0\geq \delta>0$. We obtain the global classical
solution to (\ref{1.4})-(\ref{1.6}) for $\rho_0\ge0$ and
$b<\infty$ after taking the limits $\delta\rightarrow0$. Based on
the global classical solution for the case of  $b<\infty$, where
the estimates are uniform for  $b$, we can get the solution in the
exterior domain by using similar arguments as in \cite{8}.

In the section, we denote by ``$c$"  the generic constant
depending on $a$, $\|\rho_0\|_{H^2}$, $\|\rho_0^\gamma\|_{H^2}$,
$\|{\bf u}_0\|_{D_0^1}$, $\|{\bf u}_0\|_{D^3}$, $T$ and some other
known constants but independent of $\delta$ and $b$.

Before proving Theorem \ref{th:1.1}, we give the following
auxiliary theorem.
\begin{Theorem}\label{th:2.1} Consider the same assumptions as in Theorem \ref{th:1.1},
and in addition assume that $\rho_0\geq\delta>0$ and $b<\infty$.
Then for any $T>0$,  there exists a unique global classical
solution $(\rho, u)$ to (\ref{1.4})-(\ref{1.6}) satisfying \beqno
&\rho\in C([0,T];H^2(I)), \ \ \rho\ge\frac{\D\delta}{c},\
u_{tt}\in L^2([0,T];L^2(I)), & \\&\ \ u\in C([0,T];H^3(I)\cap
H^1_0(I)), \ \ u_t\in C([0,T];H_0^1(I))\cap L^2([0,T];H^2(I)).&
\eeqno
\end{Theorem}
The local solution  in Theorem \ref{th:2.1} can be obtained by the
successive approximations as in \cite{8,10},  we omit it here for
simplicity. The regularities guarantee the uniqueness (refer for
instance to \cite{7,8}). Based on it, Theorem \ref{th:2.1} can be
proved by some {\it a priori} estimates globally in time.

For $T\in(0,+\infty)$, let $(\rho,u)$ be the classical solution of
(\ref{1.4})-(\ref{1.6}) as in Theorem \ref{th:2.1}. Then we have
the following estimates (cf. \cite{9} and \cite{11}):
\begin{Lemma} \label{le:2.1}For any $0\le t\le T$, it holds that
$$\|(\rho,\rho^\gamma)\|_{H_r^2}+\|(\rho_t,(\rho^\gamma)_t)\|_{H_r^1}+\int_0^T\|(\rho_{tt},(\rho^\gamma)_{tt})\|_{L_r^2}^2\le
c,\ \rho\ge\frac{\delta}{c},$$ and
$$
\int_I(r^m\rho
u_t^2+r^mu_{rr}^2+r^mu_r^2+r^{m-2}u^2)+\int_{Q_T}(r^mu_{rt}^2+r^{m-2}u_t^2+r^mu_{rrr}^2)\le
c,
$$
where $ \|(h_1,h_2)\|_X=\|h_1\|_X+\|h_2\|_X,$ for some Banach
space X, and $h_i\in X$, i=1,2.
\end{Lemma}

\begin{Remark}
 $(i)$ For the estimates about $\rho^\gamma$,  since $\rho$ and $\rho^\gamma$  satisfy the linear transport
equations: $\rho_t+{\bf u}\cdot \nabla \rho +\rho \nabla\cdot {\bf
u}=0$ and  $(\rho^\gamma)_t+{\bf u}\cdot \nabla (\rho^\gamma)
+\gamma\rho^\gamma \nabla\cdot {\bf u}=0$ respectively, then by
using the similar arguments as that of Lemma 3.6  in \cite{9}, we
get
\begin{eqnarray*}
&& \frac{d}{dt}\left\{\int_\Omega |\nabla^2\rho|^2+\int_\Omega
|\nabla^2(\rho^\gamma)|^2\right\}\\
&\leq & c\|\nabla {\bf u}(\cdot, t)\|_{H^1}\left(\|\nabla
\rho(\cdot,
t)\|_{H^1}^2+ \|\nabla (\rho^\gamma)(\cdot, t)\|_{H^1}^2\right)\\
&& +c\left(\|\nabla^2 G(\cdot,
t)\|_{L^2}+\|\nabla^2(\rho^\gamma)(\cdot,
t)\|_{L^2}\right)\left(\|\nabla^2\rho(\cdot,
t)\|_{L^2}+\|\nabla^2(\rho^\gamma)(\cdot, t)\|_{L^2}\right)\\
&\leq& c\left(\|\nabla \rho(\cdot, t)\|_{H^1}^2+\|\nabla
(\rho^\gamma)(\cdot, t)\|_{H^1}^2+\|G(\cdot, t)\|_{H^2}^2\right),
\end{eqnarray*}
where $G=\nu \nabla\cdot {\bf u}-\rho^\gamma$ is the effective
viscous flux; we have used the estimates $\|{\bf u}(\cdot,
t)\|_{D_0^1}\leq c$, $\|{\bf u}(\cdot, t)\|_{D^2}\leq c$ and
$\|\rho(\cdot, t)\|_{H^1}\leq c$ in \cite{9}. Since $\|G(\cdot,
t)\|_{H^2}\in L^2(0, T)$ given in \cite{9}, hence, an application
of the Gronwall inequality gives
$$\|(\rho, \rho^\gamma)(\cdot, t)\|_{H^2}\leq c,\eqno(*)$$ where we have
used the following Sobolev inequalities for radially symmetric
functions defined in $\Omega'=\{{\bf x}\in \mathbb{R}^n;|x|>
a>0\}$: $\|\nabla\rho\|_{L^\infty}\leq c \|\nabla \rho\|_{H^1}$,
$\|\nabla {\bf u}\|_{L^\infty}\leq c \|\nabla {\bf u}\|_{H^1}$.
From $(*)$ and a direct calculation, we get
$$
\|(\rho,\rho^\gamma)\|_{H_r^2}\le c.
$$

$(ii)$ For the case $a\geq 0$ and $\gamma$ different from that in
\cite{wei}, it is interesting to investigate the existence of
global spherically symmetric strong solutions or spherically
symmetric classical solutions. We leave it as a forthcoming paper.
\end{Remark}

\begin{Lemma}\label{le:2.2} For any $0\le t\le T$, it holds
$$
\int_I(r^mu_{rt}^2+r^{m-2}u_t^2)+\int_{Q_T}r^m\rho u_{tt}^2\le c.
$$
\end{Lemma}
\noindent{\it Proof.} From (\ref{1.4}), we get \be\label{2.1} \rho
u_t+\rho uu_r+(\rho^\gamma)_r=\left(u_r+\frac{mu}{r}\right)_r+\rho
f. \ee Differentiating (\ref{2.1}) with respect to $t$, we have
\be\label{2.2} \rho u_{tt}+\rho_tu_t+\rho_tuu_r+\rho u_tu_r+\rho
uu_{rt}+(\rho^\gamma)_{rt}=u_{rrt}+mr^{-2}(ru_{rt}-u_t)+\rho_tf+\rho
f_t. \ee Multiplying (\ref{2.2}) by $r^mu_{tt}$, integrating over
$I$, and using integration by parts, Lemma \ref{le:2.1} and Cauchy
inequality, we have \beqno &&\int_Ir^m\rho
u_{tt}^2+\frac{1}{2}\frac{d}{dt}\int_I(r^mu_{rt}^2+mr^{m-2}u_t^2)
\\&=&\int_Ir^m(\rho_tf+\rho
f_t)u_{tt}-\int_Ir^m[\rho_tu_t+\rho_tuu_r+\rho u_tu_r+\rho
uu_{rt}+(\rho^\gamma)_{rt}]u_{tt}\\&\leq&-\int_Ir^m\rho_tu_tu_{tt}-\int_Ir^m(\rho_tuu_r-\rho_tf)u_{tt}+c\int_Ir^m\rho
u_t^2+\frac{1}{2}\int_Ir^m\rho
u_{tt}^2+c\int_Ir^mu_{rt}^2\\&&-\int_Ir^m(\rho^\gamma)_{rt}u_{tt}+c\int_Ir^mf_t^2,
\eeqno where we have used Lemma \ref{le:2.1} and  the following
inequality found in \cite{9}: \be\label{2.3} \sup\limits_{a\le
r\le b}|u(r,t)|\le
c\left[\int_I(r^mu_r^2+mr^{m-2}u^2)(r,t)dr\right]^\frac{1}{2},\
\textrm{for}\ u(a,t)=0, \ee and the one-dimensional Sobolev
inequality: \be\label{2.4} \sup\limits_{a\le r\le
b}|\varphi(r)|\le c\|\varphi\|_{H^1(I)}. \ee
 Thus \beqno
 &&\frac{1}{2}\int_Ir^m\rho
u_{tt}^2+\frac{1}{2}\frac{d}{dt}\int_I(r^mu_{rt}^2+mr^{m-2}u_t^2)\\
&\leq&-\int_Ir^m\rho_tu_tu_{tt}-\int_Ir^m(\rho_tuu_r-\rho_tf)u_{tt}
-\int_Ir^m(\rho^\gamma)_{rt}u_{tt}
+c\int_Ir^mu_{rt}^2+c\int_Ir^mf_t^2+c\\
 &=&-\frac{d}{dt}\int_I\left[\frac{1}{2}r^m\rho_tu_t^2
+r^m(\rho_tuu_r-\rho_tf)u_t\right]
+\frac{1}{2}\int_Ir^m\rho_{tt}u_t^2\\
&&+\int_Ir^m(\rho_{tt}uu_r+\rho_tu_tu_r+\rho_tuu_{rt}-\rho_{tt}f-\rho_tf_t)u_t-\int_Ir^m(\rho^\gamma)_{rt}u_{tt}
+c\int_Ir^mu_{rt}^2\\&&+c\int_Ir^mf_t^2+c\\
&\le&-\frac{d}{dt}\int_I\left[\frac{1}{2}r^m\rho_tu_t^2
+r^m(\rho_tuu_r-\rho_tf)u_t\right]
 -\frac{1}{2}\int_I(r^m\rho
u)_{rt}u_t^2+c\|u_t\|_{L^\infty}\int_Ir^m(\rho_{tt}^2+u^2u_r^2)\\
&&+c\|u_t\|_{L^\infty}^2\int_Ir^m(\rho_t^2+u_r^2)
+c\|u_t\|_{L^\infty}\left(\int_Ir^m\rho_t^2\right)^\frac{1}{2}\left(\int_Ir^mu_{rt}^2\right)^\frac{1}{2}\\&&+\|u_t\|_{L^\infty}\int_Ir^m(\rho_{tt}^2+f^2+\rho_t^2+f_t^2)
 -\int_Ir^m(\rho^\gamma)_{rt}u_{tt}+c\int_Ir^mu_{rt}^2+c\int_Ir^mf_t^2+c,
\eeqno where we have used $(\ref{1.4})_1$, (\ref{2.3}), Cauchy
inequality, H$\ddot{{\rm o}}$lder inequality  and Lemma
\ref{le:2.1}. Using (\ref{2.3}), Cauchy inequality, Lemma
\ref{le:2.1} again, along with (\ref{2.4}) and integration by
parts, we further obtain  \beqno &&\frac{1}{2}\int_Ir^m\rho
u_{tt}^2+\frac{1}{2}\frac{d}{dt}\int_I(r^mu_{rt}^2+mr^{m-2}u_t^2)
\\
&\le&-\frac{d}{dt}\int_I\left[\frac{1}{2}r^m\rho_tu_t^2
+r^m(\rho_tuu_r-\rho_tf)u_t\right]+\int_Ir^m(\rho
u_t+\rho_tu)u_tu_{rt}\\&&+c\left[\int_I(r^mu_{rt}^2+mr^{m-2}u_t^2)\right]^\frac{1}{2}\left[\int_Ir^m(\rho_{tt}^2+f^2+f_t^2)+1\right]
+c\int_I(r^mu_{rt}^2+mr^{m-2}u_t^2)\\&&-\int_Ir^m(\rho^\gamma)_{rt}u_{tt}+c\int_Ir^mf_t^2+c\\&\le&-\frac{d}{dt}\int_I\left[\frac{1}{2}r^m\rho_tu_t^2
+r^m(\rho_tuu_r-\rho_tf)u_t\right]+c\int_Ir^m(\rho^2u_t^4+\rho_t^2u^2u_t^2)
\\&&+c\int_I(r^mu_{rt}^2+mr^{m-2}u_t^2)\left[\int_Ir^m(\rho_{tt}^2+f_t^2)+1\right]-\int_Ir^m(\rho^\gamma)_{rt}u_{tt}+c\int_Ir^m(\rho_{tt}^2+f_t^2)
+c\\
&\le&-\frac{d}{dt}\int_I\left[\frac{1}{2}r^m\rho_tu_t^2
+r^m(\rho_tuu_r-\rho_tf)u_t\right]+c\|u_t\|_{L^\infty}^2\left(\int_Ir^m\rho
u_t^2+\int_Ir^m\rho_t^2\right)\\&&+c\int_I(r^mu_{rt}^2+mr^{m-2}u_t^2)\left[\int_Ir^m(\rho_{tt}^2+f_t^2)+1\right]
-\int_Ir^m(\rho^\gamma)_{rt}u_{tt}+c\int_Ir^m(\rho_{tt}^2+f_t^2)
+c\\&\le&-\frac{d}{dt}\int_I\left[\frac{1}{2}r^m\rho_tu_t^2
+r^m(\rho_tuu_r-\rho_tf)u_t\right]
+c\int_I(r^mu_{rt}^2+mr^{m-2}u_t^2)\left[\int_Ir^m(\rho_{tt}^2+f_t^2)+1\right]
\\&&-\int_Ir^m(\rho^\gamma)_{rt}u_{tt}+c\int_Ir^m(\rho_{tt}^2+f_t^2)+c.
\eeqno Integrating the above inequality over $(0, t)$, we have
\beqno &&\frac{1}{2}\int_0^t\int_Ir^m\rho
u_{tt}^2+\frac{1}{2}\int_I(r^mu_{rt}^2+mr^{m-2}u_t^2)\\
&\le&-\int_I\left[\frac{1}{2}r^m\rho_tu_t^2
+r^m(\rho_tuu_r-\rho_tf)u_t\right]+c\int_0^t\int_I(r^mu_{rt}^2+mr^{m-2}u_t^2)\left[\int_Ir^m(\rho_{tt}^2+f_t^2)+1\right]
\\&&-\int_0^t\int_Ir^m(\rho^\gamma)_{rt}u_{tt}+c\\&\le&\int_I\frac{1}{2}(r^m\rho
u)_ru_t^2+c\|u_t\|_{L^\infty}+c\int_0^t\int_I(r^mu_{rt}^2+mr^{m-2}u_t^2)\int_Ir^m(\rho_{tt}^2+f_t^2)
\\&&-\int_0^t\int_Ir^m(\rho^\gamma)_{rt}u_{tt}+c, \eeqno where we
have used $(\ref{1.4})_1$, (\ref{1.8}), (\ref{2.3}), Lemma
\ref{le:2.1}, Cauchy inequality and the following equalities:
$$
\int_I(r^mu_{rt}^2+mr^{m-2}u_t^2)|_{t=0}=c\int_\Omega|\nabla{\bf
u}_t|^2|_{t=0},
$$
 \beqno {\bf u}_t|_{t=0}&=& {\rho_0}^{-1}\left(L{\bf u}_0+\rho_0{\bf
f}(0)-\nabla P(\rho_0)-\rho_0{\bf u}_0\cdot\nabla{\bf u}_0\right)
\\&=&{\bf g}_1-{\bf u}_0\cdot\nabla{\bf u}_0\in D_0^1.
\eeqno Using integration by parts, (\ref{2.3}), (\ref{2.4}),
H\"older inequality, Cauchy inequality and Lemma $\ref{le:2.1}$,
we obtain \beqno &&\frac{1}{2}\int_0^t\int_Ir^m\rho
u_{tt}^2+\frac{1}{2}\int_I(r^mu_{rt}^2+mr^{m-2}u_t^2)
\\&\le&-\int_Ir^m\rho u
u_tu_{rt}+c\left[\int_I(r^mu_{rt}^2+mr^{m-2}u_t^2)\right]^\frac{1}{2}+c\int_0^t\int_I(r^mu_{rt}^2+mr^{m-2}u_t^2)\int_Ir^m(\rho_{tt}^2+f_t^2)
\\&&-\int_0^t\int_Ir^m(\rho^\gamma)_{rt}u_{tt}+c\\&\le&\frac{1}{4}\int_I(r^mu_{rt}^2+mr^{m-2}u_t^2)+c\int_0^t\int_I(r^mu_{rt}^2+mr^{m-2}u_t^2)\int_Ir^m(\rho_{tt}^2+f_t^2)
-\int_0^t\int_Ir^m(\rho^\gamma)_{rt}u_{tt}+c.\eeqno Therefore,
\beqno &&\int_0^t\int_Ir^m\rho
u_{tt}^2+\frac{1}{2}\int_I(r^mu_{rt}^2+mr^{m-2}u_t^2)\\&\le&c\int_0^t\int_I(r^mu_{rt}^2+mr^{m-2}u_t^2)\int_Ir^m(\rho_{tt}^2+f_t^2)
-\int_0^t\int_Ir^m(\rho^\gamma)_{rt}u_{tt}+c.\ \ \ \ \ \ \ \ \ \ \
\ \ \ \ \ \ \ \ \ (*_1)
 \eeqno
We shall handle the second integral of the right hand side of the
above inequality. Using integration by parts twice formally
(firstly with respect to $t$ and then with respect to $r$), we
have \beqno
-\int_0^t\int_Ir^m(\rho^\gamma)_{rt}u_{tt}&=&-\int_Ir^m(\rho^\gamma)_{rt}u_{t}+\int_Ir^m{(\rho^\gamma)_{rt}u_{t}}\big|_{t=0}
+\int_0^t\int_Ir^m(\rho^\gamma)_{rtt}u_{t}\\&\le&\int_Ir^m(\rho^\gamma)_{t}u_{rt}+m\int_Ir^{m-1}(\rho^\gamma)_{t}u_{t}
-\int_0^t\int_Ir^m(\rho^\gamma)_{tt}u_{rt}\\&&-m\int_0^t\int_Ir^{m-1}(\rho^\gamma)_{tt}u_{t}+c\\&\le&\frac{1}{4}\int_I(r^mu_{rt}^2+mr^{m-2}u_t^2)
+c. \ \ \ \ \ \ \ \ \ \ \ \ \ \ \ \ \ \ \ \ \ \ \ \ \ \ \ \ \ \ \
\ \ \ \ \ \ \ \ \ \ \ \ \ \ (*_2)\eeqno Rigorously, this process
can be done by the usual method (mollification+taking limits).
More precisely,
 we can take
place of $(\rho^\gamma)_t$ by
$J_\epsilon*\overline{(\rho^\gamma)_t}$, and then take limits
$\epsilon\rightarrow 0^+$, where
$J_\epsilon(\cdot)=\frac{1}{\epsilon}J(\frac{\cdot}{\epsilon})$;
$J(\cdot)$ is the usual mollifier in $\mathbb{R}^1$;
$\overline{(\rho^\gamma)_t}$ is an extension of $(\rho^\gamma)_t$
w.r.t. $t$, defined as follows
$$
\overline{(\rho^\gamma)_t}=\left\{\begin{array}{l}
(\rho^\gamma)_t(0), \ \ \ {\rm if} \ \ \ t\in[-1,0],
\\ [3mm](\rho^\gamma)_t(t), \ \ \ {\rm if}\ \ \ t\in[0,T],
\\ [3mm] (\rho^\gamma)_t(T), \ \ \ {\rm if}\ \ \ t\in[T, T+1].
\end{array}
\right.
$$
Substituting ($*_2$) into ($*_1$), and using Gronwall inequality,
we get
 \beqno\int_{Q_T}r^m\rho
u_{tt}^2+\int_I(r^mu_{rt}^2+r^{m-2}u_t^2)\le c. \eeqno The proof
of
  Lemma \ref{le:2.2} is complete. \qed
\begin{Lemma}\label{le:2.3} For any $0\le t\le T$, it holds
$$
\int_Ir^mu_{rrr}^2\le c.
$$
\end{Lemma}
\noindent{\it Proof.} Differentiating (\ref{2.1}) with respect to
$r$, we get \beq\label{2.5}\nonumber u_{rrr}=\rho_ru_t+\rho
u_{rt}+\rho_ruu_r+\rho u_r^2+\rho uu_{rr}+(\rho^\gamma)_{rr}
-\frac{mu_{rr}}{r}+\frac{2m(ru_r-u)}{r^3}-\rho_rf-\rho f_r.\\ \eeq
By (\ref{2.5}), (\ref{2.3}), (\ref{2.4}) and Lemma \ref{le:2.1},
we obtain \beqno\int_Ir^mu_{rrr}^2&\le&
c\int_I(r^mu_{rt}^2+mr^{m-2}u_t^2)+c\int_Ir^m|(\rho^\gamma)_{rr}|^2+c\int_Ir^m(f^2+f_r^2)+c\\&\le
&c. \eeqno The proof of   Lemma \ref{le:2.3} is complete.\qed
\begin{Lemma}\label{le:2.4} For any $0\le t\le T$, it holds
$$
\int_Ir^m(\rho_{tt}^2+|(\rho^\gamma)_{tt}|^2)+\int_{Q_T}r^mu_{rrt}^2\le
c.
$$
\end{Lemma}
\noindent{\it Proof.} By (\ref{2.2}), (\ref{2.3}), (\ref{2.4}),
Lemma \ref{le:2.1}  and Lemma \ref{le:2.2}, we get
\beqno\int_{Q_T}r^mu_{rrt}^2&\le&c\int_{Q_T}r^m\rho_t^2u_t^2+c\int_{Q_T}r^m\rho^2u_{tt}^2+c\int_{Q_T}r^m\rho_t^2u^2u_r^2+c\int_{Q_T}r^m
\rho^2u_t^2u_r^2\\&&+c\int_{Q_T}r^m\rho^2u^2u_{rt}^2+c\int_{Q_T}r^m|(\rho^\gamma)_{rt}|^2+c\int_{Q_T}(r^mu_{rt}^2+mr^{m-2}u_t^2)
\\&&+c\int_{Q_T}r^m(\rho_t^2f^2+\rho^2f_t^2)\\&\le&c. \eeqno
Multiplying $(\ref{1.4})_1$ by $\gamma\rho^{\gamma-1}$,we get
\be\label{2.6} (\rho^\gamma)_t+\gamma\rho^\gamma
u_r+(\rho^\gamma)_ru+m\gamma r^{-1}\rho^\gamma u=0. \ee
$(\ref{1.4})_1$, (\ref{2.3}), (\ref{2.4}), (\ref{2.6}), Lemma
\ref{le:2.1}  and Lemma \ref{le:2.2} imply
$$
\int_Ir^m(\rho_{tt}^2+|(\rho^\gamma)_{tt}|^2)\le c.
$$
This proves Lemma \ref{le:2.4}. \qed

To sum up, we get
\be\label{2.7}\|(\rho,\rho^\gamma)\|_{H_r^2}+\|(\rho_t,(\rho^\gamma)_t)\|_{H_r^1}+\|(\rho_{tt},(\rho^\gamma)_{tt})\|_{L_r^2}\le
c,\ \rho\ge\frac{\delta}{c}, \ee and

 \be\label{2.8} \int_Ir^m(u_{rt}^2+r^{-2}u_t^2+u_{rrr}^2+u_{rr}^2+u_r^2+r^{-2}u^2)+\int_{Q_T}r^m(\rho
u_{tt}^2+u_{rrt}^2)\le c, \ee By (\ref{2.7}) and (\ref{2.8}), we
complete the proof of Theorem \ref{th:2.1}.\qed

\vspace{4mm}

{\bf Proof of Theorem \ref{th:1.1}}:

Let $b<\infty$, and denote $\rho_0^\delta=\rho_0+\delta$, for
$\delta\in (0,1)$, we have \be\label{2.9}
\rho_0^\delta\rightarrow\rho_0, \ \ {\rm in} \ \ H^2(I), \ee
\be\label{2.10} (\rho_0^\delta)^\gamma\rightarrow\rho_0^\gamma, \
\ {\rm in} \ \ H^2(I). \ee Let $u_0^\delta$ be the unique solution
to the equation: \be\label{2.11}
\left(u_{0r}^\delta+\frac{mu_0^\delta}{r}\right)_r-[(\rho_0^\delta)^\gamma]_r=\rho_0^\delta[g_1-f(
0)],\ \textrm{in}\ I, \ee for $u_0^\delta|_{\partial I}=0$.
(\ref{1.8}) implies
 \be\label{2.12}
(u_{0r}+\frac{mu_0}{r})_r-[(\rho_0)^\gamma]_r=\rho_0[g_1-f(0)],\
\textrm{in}\ I, \ee for $u_0|_{\partial I}=0$.

\vspace{2mm}

From (\ref{2.9})-(\ref{2.12}) and the standard elliptic estimates,
we obtain \be\label{2.13} u_0^\delta\rightarrow u_0, \ \ {\rm in}
\ \ H^3(I), \ \ {\rm as} \ \ \delta\rightarrow0.\ee

 Consider (\ref{1.4})-(\ref{1.6}) with initial-boundary data
 replaced by
$$
(\rho^\delta,u^\delta)|_{t=0}=(\rho_0^\delta,u_0^\delta), \ \ {\rm
in} \ \ I,
$$
and
$$
u^\delta|_{\partial I}=0, \ \ {\rm for} \ \ t\ge0.
$$
Then we get a unique solution $(\rho^\delta,u^\delta)$ for each
$\delta>0$ by Theorem \ref{th:2.1}.

Following the estimates in the proofs of Theorem \ref{th:2.1}, we
can also get (\ref{2.7}) and (\ref{2.8}) with $(\rho,u)$ replaced
by $(\rho^\delta,u^\delta)$. Since $b<\infty$ and $a>0$, it
follows from (\ref{2.7}) and (\ref{2.8}) that
\be\label{2.14}\|(\rho^\delta,(\rho^\delta)^\gamma)\|_{H^2(I)}+\|(\rho^\delta_t,((\rho^\delta)^\gamma)_t)\|_{H^1(I)}
+\|(\rho^\delta_{tt},((\rho^\delta)^\gamma)_{tt})\|_{L^2(I)}\le
c,\ee and \be\label{2.15} \rho^\delta\ge\frac{\delta}{c},\
\|u^\delta_t\|_{H^1(I)}+\|u^\delta\|_{H^3(I)}+\int_{Q_T}(\rho^\delta
|u_{tt}^\delta|^2+|u_{rrt}^\delta|^2)\le c(b). \ee Based on the
estimates  in (\ref{2.14}) and (\ref{2.15}), we get a solution
$(\rho,u)$ to (\ref{1.4})-(\ref{1.6}) after taking the limit
$\delta\rightarrow0$ (take the subsequence if necessary),
satisfying \beq\label{2.16}\begin{cases} (\rho,\rho^\gamma)\in
L^\infty([0,T];H^2(I)), \ \ \ (\rho_t,(\rho^\gamma)_t)\in
L^\infty([0,T];H^1(I)), \ \ \
\\ (\rho_{tt},(\rho^\gamma)_{tt})\in L^\infty([0,T];L^2(I)), \ \ \
\rho\geq0,\ \ (\rho u_t)_t\in L^2([0,T];L^2(I)),  \\
u\in L^\infty([0,T];H^3(I)\cap H^1_0(I)), \ \ \ u_t\in
L^\infty([0,T];H_0^1(I))\cap L^2([0,T];H^2(I)).
\end{cases} \eeq
Since $u\in L^\infty([0,T];H^3(I))$ and $u_t\in
L^\infty([0,T];H_0^1(I))$, then we get  $ u\in C([0,T];H^2(I))$
(refer to \cite{25}). By $(\ref{1.4})_1$, (\ref{2.6}) and similar
arguments as in \cite{9, 10}, we get \be\label{2.17} \rho\in
C([0,T];H^2(I)),\ \rho_t\in C([0,T];H^1(I)), \ee and
 \be\label{2.18} \rho^\gamma\in
C([0,T];H^2(I)),\ (\rho^\gamma)_t\in C([0,T];H^1(I)). \ee Denote
$G=(u_r+\frac{mu}{r}-\rho^\gamma)_r+\rho f$. By (\ref{2.1}) and
(\ref{2.16}), we have
$$ G=\rho u_t+\rho uu_r\in
L^2([0,T];H^2(I)),$$
$$ G_t=(\rho u_t+\rho uu_r)_t\in L^2([0,T];L^2(I)).
$$
By the embedding theorem (\cite{25}), we have $G\in
C([0,T];H^1(I))$. Since $\rho f\in C([0,T];H^1(I))$, we get
$$
\left(u_r+\frac{mu}{r}-\rho^\gamma\right)_r\in C([0,T];H^1(I)).
$$
This means \be\label{2.19} u_r+\frac{mu}{r}-\rho^\gamma\in
C([0,T];H^2(I)). \ee By (\ref{2.18}) and (\ref{2.19}), we get
$$
u_r+\frac{mu}{r}\in C([0,T];H^2).
$$
This together with $ u\in C([0,T];H^2(I))$ implies \be\label{2.20}
u\in C([0,T];H^3(I)). \ee $(\ref{1.4})_2$, (\ref{2.17}),
(\ref{2.18}) and (\ref{2.20}) give \be\label{2.21} (\rho u)_t\in
C([0,T];H^1(I)). \ee Denote \be\label{2.22}\rho({\bf
x},t)=\rho(r,t),\ {\bf u}({\bf x},t)=u(r,t)\frac{{\bf x}}{r}.\ee
It follows from (\ref{2.16})-(\ref{2.18}) and
(\ref{2.20})-(\ref{2.22}), we complete the proof of Theorem
\ref{th:1.1} for $b<\infty$. For $b=\infty$, we can use the
similar methods as in \cite{8} together with the estimates
(\ref{2.7}) and (\ref{2.8}) uniform for $b$ to get it. We omit
details here for simplicity. \qed

 \setcounter{section}{2} \setcounter{equation}{0}
\section{ \ Proof of Theorem \ref{th:1.2}}

Similar to the proof of Theorem \ref{th:1.1}, we need the
following auxiliary theorem.
\begin{Theorem}\label{th:3.1} Consider the same assumptions as in Theorem \ref{th:1.2},
and in addition assume that $\rho_0\geq\delta>0$ and $b<\infty$.
Then for any $T>0$, there exists a unique global classical
solution $(\rho, u)$ to (\ref{1.4})-(\ref{1.6}) satisfying \beqno
&\rho\in C([0,T];H^5(I)),\ \rho\ge\frac{\D\delta}{c}, \ u\in
C([0,T];H^5(I)\cap H^1_0(I))\cap L^2([0,T];H^6(I)),&\\&u_t\in
C([0,T];H^3(I)\cap H^1_0(I))\cap L^2([0,T];H^4(I)), &\\    \ & \
u_{tt}\in C([0,T];H_0^1(I))\cap L^2([0,T];H^2(I)), \ u_{ttt}\in
 L^2([0,T];L^2(I)).& \eeqno
\end{Theorem}
 The proof of local existence of the solutions
as in Theorem \ref{th:3.1} can be done by the successive
approximations as in \cite{10} and references therein, together
with the estimates in Section 3 for higher order derivatives of
the solutions. We omit it here for brevity. Therefore, Theorem
\ref{th:3.1} can be proved by some {\it a priori} estimates
globally in time. Since (\ref{2.7}) and (\ref{2.8}) are also valid
here, we need some other {\it a priori} estimates about higher
order derivatives of $(\rho,u)$. The generic positive constant $c$
may depend on the initial data presented in Theorem \ref{th:1.2}
and other known constants but independent of $\delta$ and $b$.

\begin{Lemma}\label{le:3.1} For any $0\le t\le T$, it holds
$$
\int_Ir^m[\rho_{rrr}^2+\rho_{rrt}^2+|(\rho^\gamma)_{rrr}|^2
+|(\rho^\gamma)_{rrt}|^2]+\int_{Q_T}r^m[\rho_{rtt}^2+|(\rho^\gamma)_{rtt}|^2+u_{rrrr}^2]\le
c.
$$
\end{Lemma}
\noindent{\it Proof.}
 Taking derivative of order three on both sides of $(\ref{1.4})_1$
with respect to $r$, multiplying it by $r^m\rho_{rrr}$, and
integrating by parts over $I$, we have \beqno
&&\frac{1}{2}\frac{d}{dt}\int_Ir^m\rho_{rrr}^2
=\int_Ir^m\rho_{rrr}\Big[-mr^{-1}\rho
u_{rrr}-3mr^{-1}\rho_ru_{rr}-3mr^{-1}\rho_{rr}u_r+\frac{3m\rho
u_{rr}}{r^2}\\&&-\frac{6m\rho
u_r}{r^3}-mr^{-1}\rho_{rrr}u+\frac{3m\rho_{rr}
u}{r^2}+\frac{6m\rho_r u_r}{r^2}-\frac{6m\rho_r
u}{r^3}+\frac{6m\rho u}{r^4}-4\rho_{rrr}u_r-6\rho_{rr}u_{rr}\\&&\
\ \ \ \ \ \ \ \ \ -4\rho_ru_{rrr}\Big]+\frac{1}{2}\int_Ir^m
\rho_{rrr}^2
 u_r+\frac{1}{2}\int_Imr^{m-1} \rho_{rrr}^2 u-\int_Ir^m\rho\rho_{rrr}
u_{rrrr}. \eeqno By Sobolev inequality, (\ref{2.3}), (\ref{2.4})
(\ref{2.7}), (\ref{2.8})  and Cauchy inequality, we get
\be\label{3.2} \frac{d}{dt}\int_Ir^m\rho_{rrr}^2\le
c\int_Ir^m\rho_{rrr}^2+c\int_Ir^mu_{rrrr}^2+c.\ee Similarly to
(\ref{3.2}), we get from (\ref{2.6}) \be\label{3.3}
\frac{d}{dt}\int_Ir^m|(\rho^\gamma)_{rrr}|^2\le c\int_I
r^m|(\rho^\gamma)_{rrr}|^2+c\int_Ir^mu_{rrrr}^2+c. \ee By
(\ref{3.2}) and (\ref{3.3}), we have \be\label{3.4}
\frac{d}{dt}\int_Ir^m(\rho_{rrr}^2+|(\rho^\gamma)_{rrr}|^2)\le
c\int_Ir^m(\rho_{rrr}^2+|(\rho^\gamma)_{rrr}|^2)+c\int_Ir^mu_{rrrr}^2+c.
\ee Differentiating (\ref{2.5}) with respect to $r$, we have
\beq\label{3.5} \nonumber u_{rrrr}&=&\rho_{rr}u_t+2\rho_ru_{rt}
+\rho u_{rrt}+(\rho_ruu_r+\rho u_r^2+\rho
uu_{rr})_r+(\rho^\gamma)_{rrr} -mr^{-1}u_{rrr}\\&&
+3mr^{-2}u_{rr}-6mr^{-3}u_r+\frac{6mu}{r^4}-\rho_{rr}f-2\rho_rf_r-\rho
f_{rr}. \eeq From (\ref{3.5}), (\ref{2.3}), (\ref{2.4}),
(\ref{2.7}) and (\ref{2.8}), we obtain \be\label{3.6}
\int_Ir^mu_{rrrr}^2 \le
c\int_Ir^m|(\rho^\gamma)_{rrr}|^2+c\int_Ir^m\rho
u_{rrt}^2+c\int_Ir^mf_{rr}^2+c. \ee By (\ref{3.4}), (\ref{3.6}),
(\ref{2.4}) and (\ref{2.7}), we get \beqno
\frac{d}{dt}\int_Ir^m(\rho_{rrr}^2+|(\rho^\gamma)_{rrr}|^2)\le
c\int_Ir^m(\rho_{rrr}^2+|(\rho^\gamma)_{rrr}|^2)+c\int_Ir^mu_{rrt}^2+c\int_Ir^mf_{rr}^2+c.
 \eeqno
By Gronwall inequality  and (\ref{2.8}), we obtain \be\label{3.7}
\int_Ir^m(\rho_{rrr}^2+|(\rho^\gamma)_{rrr}|^2)\le c. \ee It
follows from $(\ref{1.4})_1$, (\ref{2.3}), (\ref{2.4}),
(\ref{2.6}), (\ref{2.7}), (\ref{2.8}), (\ref{3.6}) and (\ref{3.7})
that
$$
\int_Ir^m[\rho_{rrt}^2
+|(\rho^\gamma)_{rrt}|^2]+\int_{Q_T}r^m[\rho_{rtt}^2+|(\rho^\gamma)_{rtt}|^2+u_{rrrr}^2]\le
c.
$$
The proof of Lemma \ref{le:3.1} is complete. \qed
\begin{Lemma}\label{le:3.2} For any $T>0$, we have
$$
\|(\sqrt{\rho})_r\|_{L^\infty(Q_T)}+\|(\sqrt{\rho})_t\|_{L^\infty(Q_T)}\le
c.
$$
\end{Lemma}
\noindent{\it Proof.} Multiplying $(\ref{1.4})_1$ by
$\D\frac{1}{2\sqrt{\rho}}$, we have \be\label{3.8}
(\sqrt{\rho})_t+\frac{mr^{-1}}{2}\sqrt{\rho}u+\frac{1}{2}\sqrt{\rho}u_r+(\sqrt{\rho})_ru=0.
\ee Differentiating (\ref{3.8}) with respect to $r$, we get
$$
(\sqrt{\rho})_{rt}+\frac{mr^{-1}}{2}(\sqrt{\rho})_ru+\frac{mr^{-1}}{2}\sqrt{\rho}u_r-\frac{m\sqrt{\rho}u}{2r^2}+\frac{3}{2}(\sqrt{\rho})_ru_r+
\frac{1}{2}\sqrt{\rho}u_{rr}+(\sqrt{\rho})_{rr}u=0.
$$
Denote $h=(\sqrt{\rho})_r$, we have
$$
h_t+h_ru+h(\frac{mr^{-1}}{2}u+\frac{3}{2}u_r)+\frac{mr^{-1}}{2}\sqrt{\rho}u_r-\frac{m\sqrt{\rho}u}{2r^2}+
\frac{1}{2}\sqrt{\rho}u_{rr}=0,
$$
which implies \beq\label{3.9}
&&\nonumber\frac{d}{dt}\left\{h\exp\left[\int_0^t(\frac{mr^{-1}}{2}u+
\frac{3}{2}u_r)(r(\tau,y),\tau)d\tau\right]\right\}\\&=&-\left(\frac{mr^{-1}\sqrt{\rho}u_r}{2}-\frac{m\sqrt{\rho}u}{2r^2}+
\frac{\sqrt{\rho}u_{rr}}{2}\right)\nonumber\\
&&\ \ \ \ \ \times \exp\left[\int_0^t\left(\frac{mr^{-1}u}{2}+
\frac{3u_r}{2}\right)(r(\tau,y),\tau)d\tau\right], \eeq where
$r(t,y)$ satisfies
 $$
\begin{cases}
\frac{d r(t,y)}{d t}=u(r(t,y),t),\ 0\le t<s,\\
r(s,y)=y.
\end{cases} $$
Integrating (\ref{3.9}) over $(0,s)$, we get
$$
\arraycolsep=1.5pt
\begin{array}{rl}
h(y,s)= & \displaystyle
\exp\left(-\int_0^s\left(\frac{mr^{-1}u}{2}+
\frac{3u_r}{2}\right)(r(\tau,y),\tau)d\tau\right)h(r(0,y),0)
\\ [5mm]
& \displaystyle
-\int_0^s\Big[(\frac{mr^{-1}\sqrt{\rho}u_r}{2}-\frac{m\sqrt{\rho}u}{2r^2}+
\frac{\sqrt{\rho}u_{rr}}{2})\\
&\ \ \ \  \ \ \ \ \ \
\D\times\exp\left(\int_s^t\left(\frac{mr^{-1}u}{2}+
\frac{3u_r}{2}\right)(r(\tau,y),\tau)d\tau\right)\Big]dt.
\end{array}
$$
This together with (\ref{2.3}), (\ref{2.4}), (\ref{2.7}) and
(\ref{2.8}) implies \be\label{3.10}
\|(\sqrt{\rho})_r\|_{L^\infty(Q_T)}\le c. \ee From (\ref{3.8}),
(\ref{3.10}), (\ref{2.3}), (\ref{2.4}), (\ref{2.7}) and
(\ref{2.8}), we get
$$
\|(\sqrt{\rho})_t\|_{L^\infty(Q_T)}\le c.
$$
The proof of Lemma \ref{le:3.2} is complete. \qed

\vspace{3mm}

\begin{Lemma}\label{le:3.3}For any $0\le t\le T$, it holds
$$
\int_Ir^m(\rho^3 u_{tt}^2+\rho
u_{rrt}^2+u_{rrrr}^2)+\int_{Q_T}r^m(\rho^2
u_{rtt}^2+u_{rrrt}^2)\le c.
$$
\end{Lemma}
\noindent{\it Proof.} Differentiating (\ref{2.2}) with respect to
$t$, multiplying it by $r^m\rho^2u_{tt}$, and integrating by parts
over $I$, we have \beqno &&\frac{1}{2}\frac{d}{dt}\int_Ir^m\rho^3
u_{tt}^2+\int_I(r^m\rho^2u_{rtt}^2+mr^{m-2}\rho^2u_{tt}^2)\\&=&-\frac{1}{2}\int_Ir^m\rho^2\rho_tu_{tt}^2-\int_Ir^m\rho^2u_{tt}\Big[\rho_{tt}u_t
+\rho_{tt}uu_r+2\rho_tu_tu_r+2\rho_tuu_{rt}+\rho
u_{tt}u_r\\&&+2\rho u_tu_{rt}+\rho
uu_{rtt}+(\rho^\gamma)_{rtt}\Big]-2\int_Ir^m\rho\rho_ru_{tt}u_{rtt}
 +\int_Ir^m\rho^2u_{tt}(\rho_{tt}f+2\rho_tf_t+\rho
f_{tt})
\\&\le&c\int_Ir^m\rho
u_{tt}^2+\frac{1}{4}\int_Ir^m\rho^2u_{rtt}^2+c\int_Ir^m|(\rho^\gamma)_{rtt}|^2
 -4\int_Ir^m\rho
u_{rtt}\sqrt{\rho}u_{tt}(\sqrt{\rho})_r\\&&+c\int_Ir^m(\rho_{tt}^2f^2+\rho_t^2f_t^2+\rho^2f_{tt}^2)+c
\\&\le&c\int_Ir^m\rho
u_{tt}^2+\frac{1}{4}\int_Ir^m\rho^2u_{rtt}^2+c\int_Ir^m|(\rho^\gamma)_{rtt}|^2
+\frac{1}{4}\int_Ir^m\rho^2u_{rtt}^2+\int_Ir^mf_{tt}^2+c, \eeqno
where we have used (\ref{2.3}), (\ref{2.4}), (\ref{2.7}),
(\ref{2.8}), Lemma \ref{le:3.2} and Cauchy inequality.

Thus, \beq\label{3.12} &&\nonumber\frac{d}{dt}\int_Ir^m\rho^3
u_{tt}^2+\int_I(r^m\rho^2u_{rtt}^2+mr^{m-2}\rho^2u_{tt}^2)
\\&\le&c\int_Ir^m\rho
u_{tt}^2+c\int_Ir^m|(\rho^\gamma)_{rtt}|^2 +\int_Ir^mf_{tt}^2+c.
\eeq Integrating (\ref{3.12}) over $(0,t)$, and using (\ref{2.8})
and Lemma \ref{le:3.1}, we get \be\label{3.13} \int_Ir^m\rho^3
u_{tt}^2(t)+\int_0^t\int_I(r^m\rho^2u_{rtt}^2+mr^{m-2}\rho^2u_{tt}^2)\le
\int_Ir^m\rho^3 u_{tt}^2(0)+c. \ee By (\ref{2.1}), (\ref{2.2}),
(\ref{1.8}) and $\sqrt{\rho_0}\nabla^2{\bf g}_1\in L^2$, we have
\be\label{3.14} \int_Ir^m\rho^3 u_{tt}^2(0)\le c.\ee (\ref{3.13})
and (\ref{3.14}) give
 \be\label{3.15} \int_Ir^m\rho^3
u_{tt}^2+\int_{Q_T}(r^m\rho^2u_{rtt}^2+r^{m-2}\rho^2u_{tt}^2)\le
c.\ee By (\ref{2.2}), Cauchy inequality, (\ref{2.3}), (\ref{2.4}),
(\ref{2.7}), (\ref{2.8}) and (\ref{3.15}), we have
\beq\label{3.16}\nonumber \int_Ir^m\rho
u_{rrt}^2&\le&c\int_Ir^m\rho\rho_t^2u_t^2+c\int_Ir^m\rho^3u_{tt}^2+c\int_Ir^m\rho\rho_t^2u^2u_r^2
+c\int_Ir^m\rho^3u_t^2u_r^2 \\&&\nonumber +
c\int_Ir^m\rho^3u^2u_{rt}^2+c\int_Ir^m\rho|(\rho^\gamma)_{rt}|^2+c\int_Ir^{m-2}\rho
u_{rt}^2+c\int_Ir^{m-4}\rho
u_t^2\\&&\nonumber+c\int_Ir^m\rho\rho_t^2f^2+c\int_Ir^m\rho^3f_t^2\\&\le&c.
\eeq
 Differentiating (\ref{2.2}) with respect to $r$, we get
\beq\label{3.17}\nonumber u_{rrrt}&=&
\rho_{rt}u_t+2(\sqrt{\rho})_r\sqrt{\rho}u_{tt}+\rho_tu_{rt}+\rho
u_{rtt}+\rho_{rt}uu_r
+\rho_ru_tu_r+\rho_ruu_{rt}+\rho_tu_r^2\\&&\nonumber+2\rho
u_ru_{rt}+\rho_tuu_{rr}+\rho u_tu_{rr}
 +\rho
uu_{rrt}+(\rho^\gamma)_{rrt}-mr^{-1}u_{rrt}+2mr^{-2}u_{rt}\\&&-\frac{2m
u_t}{r^3}-\rho_{rt}f-\rho_tf_r-\rho_rf_t-\rho f_{rt}. \eeq
 By (\ref{3.15}), (\ref{3.17}), (\ref{2.3}), (\ref{2.4}), (\ref{2.7}), (\ref{2.8}),
 Lemma \ref{le:3.1} and Lemma \ref{le:3.2}, we have
$$
\int_{Q_T}r^m u_{rrrt}^2\le c.
$$
(\ref{3.6}), (\ref{3.16}) and Lemma \ref{le:3.1} immediately give
$$
\int_Ir^mu_{rrrr}^2\le c.
$$
 The proof of Lemma
\ref{le:3.3} is complete. \qed
\begin{Lemma}\label{le:3.4}For any $0\le t\le T$, it holds
$$
\int_Ir^m[\rho_{rrrr}^2+|(\rho^\gamma)_{rrrr}|^2]+\int_{Q_T}r^mu_{rrrrr}^2\le
c.
$$
\end{Lemma}
\noindent{\it Proof.} Differentiating (\ref{1.4})$_1$ four times
with respect to $r$, we get \beq\label{3.18}
\nonumber&&\rho_{rrrrt}+mr^{-1}\rho
u_{rrrr}+4mr^{-1}\rho_ru_{rrr}+6mr^{-1}\rho_{rr}u_{rr}-\frac{4m\rho
u_{rrr}}{r^2}+\frac{12m\rho u_{rr}}{r^3}-\frac{18m\rho
u_r}{r^4}\\&&\nonumber+4mr^{-1}\rho_{rrr}u_r+mr^{-1}\rho_{rrrr}u-\frac{4m\rho_{rrr}
u}{r^2}-\frac{12m\rho_r u_{rr}}{r^2}-\frac{12m\rho_{rr}
u_r}{r^2}+\frac{24m\rho_ru_r}{r^3}\\&&\nonumber+\frac{12m\rho_{rr}u}{r^3}-\frac{6m\rho
u_r}{r^4}-\frac{24m\rho_r u}{r^4}+\frac{24m\rho u}{r^5}+\rho
u_{rrrrr}+5\rho_ru_{rrrr}\\&&+10\rho_{rr}u_{rrr}+10\rho_{rrr}u_{rr}+5\rho_{rrrr}u_r+\rho_{rrrrr}u=0
 \eeq
Multiplying (\ref{3.18}) by $r^m\rho_{rrrr}$, integrating over
$I$, and using integration by parts, (\ref{2.3}), (\ref{2.4}),
(\ref{2.7}), (\ref{2.8}), Lemma \ref{le:3.1}, Lemma \ref{le:3.3},
we get \be\label{3.19} \frac{d}{dt}\int_Ir^m\rho_{rrrr}^2 \le
c\int_Ir^m\rho_{rrrr}^2+c\int_Ir^mu_{rrrrr}^2+c.\ee Similarly, we
have \be\label{3.20}
\frac{d}{dt}\int_Ir^m|(\rho^\gamma)_{rrrr}|^2\le
c\int_Ir^m|(\rho^\gamma)_{rrrr}|^2+c\int_Ir^mu_{rrrrr}^2+c. \ee By
(\ref{3.19}) and (\ref{3.20}), we obtain \be\label{3.21}
\frac{d}{dt}\int_Ir^m(\rho_{rrrr}^2+|(\rho^\gamma)_{rrrr}|^2)\le
c\int_Ir^m(\rho_{rrrr}^2+|(\rho^\gamma)_{rrrr}|^2)+c\int_Ir^mu_{rrrrr}^2+c.
\ee Now we estimate $\int_Ir^mu_{rrrrr}^2$. Differentiating
(\ref{3.5}) with respect to $r$, we have \beq\label{3.22}\nonumber
u_{rrrrr}&=&\rho_{rrr}u_t+3\rho_{rr}u_{rt}+3\rho_ru_{rrt}+\rho
u_{rrrt}+(\rho_ruu_r+\rho u_r^2+\rho
uu_{rr})_{rr}+(\rho^\gamma)_{rrrr}\\&&\nonumber-mr^{-1}u_{rrrr}+4mr^{-2}u_{rrr}-12mr^{-3}u_{rr}+24mr^{-4}u_r-\frac{24m
u}{r^5}
\\&& -\rho_{rrr}f-3\rho_{rr}f_r-3\rho_rf_{rr}-\rho
f_{rrr}. \eeq It follows from (\ref{3.22}), (\ref{2.3}),
(\ref{2.4}), (\ref{2.7}), (\ref{2.8}), Lemma \ref{le:3.1} and
Lemma \ref{le:3.3} that \beq\label{3.23}
\nonumber\int_Ir^mu_{rrrrr}^2&\le&
c\int_Ir^m\rho_r^2u_{rrt}^2+c\int_Ir^m\rho^2u_{rrrt}^2
+c\int_Ir^m|(\rho^\gamma)_{rrrr}|^2\\&&+c\int_Ir^m(f_r^2+f_{rr}^2+f_{rrr}^2)+c.
\eeq Since
$$
\int_Ir^m\rho_r^2u_{rrt}^2=4\int_Ir^m\rho|(\sqrt{\rho})_r|^2u_{rrt}^2.
$$
This together with (\ref{2.7}), (\ref{2.8}), Lemma \ref{le:3.2}
and Lemma \ref{le:3.3} gives  \be\label{3.24}
\int_Ir^mu_{rrrrr}^2\le c\int_Ir^m\rho^2 u_{rrrt}^2
+c\int_Ir^m|(\rho^\gamma)_{rrrr}|^2+c\int_Ir^m(f_r^2+f_{rr}^2+f_{rrr}^2)+c
\ee Substituting (\ref{3.24}) into  (\ref{3.21}), we obtain \beqno
&&\frac{d}{dt}\int_Ir^m(\rho_{rrrr}^2+|(\rho^\gamma)_{rrrr}|^2)
\\&\le& c\int_Ir^m(\rho_{rrrr}^2+|(\rho^\gamma)_{rrrr}|^2)
+c\int_Ir^mu_{rrrt}^2 +c\int_Ir^m(f_r^2+f_{rr}^2+f_{rrr}^2)+c.
\eeqno Using Gronwall inequality and Lemma \ref{le:3.3}, we get
\be\label{3.25}
\int_Ir^m(\rho_{rrrr}^2+|(\rho^\gamma)_{rrrr}|^2)\le c. \ee It
follows from (\ref{3.24}), (\ref{3.25}) and Lemma \ref{le:3.3}
 that
 $$ \int_{Q_T}r^mu_{rrrrr}^2\le
c.
$$ This proves Lemma \ref{le:3.4}. \qed

From $(\ref{1.4})_1$, (\ref{2.3}), (\ref{2.4}), (\ref{2.6}),
(\ref{2.7}), (\ref{2.8})  and Lemmas \ref{le:3.1}-\ref{le:3.4}, we
immediately get the following estimate.
\begin{Lemma}\label{le:3.5}For any $0\le t\le T$, it holds
\begin{eqnarray*}
&&\int_Ir^m\left[\rho_{rtt}^2+|(\rho^\gamma)_{rtt}|^2+\rho_{rrrt}^2
+|(\rho^\gamma)_{rrrt}|^2\right]\\
&& \ \ \ \ \
+\int_{Q_T}r^m\left[\rho_{ttt}^2+|(\rho^\gamma)_{ttt}|^2+\rho_{rrtt}^2
+|(\rho^\gamma)_{rrtt}|^2\right] \le  c.
\end{eqnarray*}
\end{Lemma}

\begin{Lemma}\label{le:3.6}For any $0\le t\le T$, it holds
$$
\int_I(r^m\rho^4u_{rtt}^2+r^{m-2}\rho^4u_{tt}^2)+\int_{Q_T}r^m\rho^5u_{ttt}^2
\le c.
$$
\end{Lemma}
\noindent{\it Proof.}  Differentiating (\ref{2.2}) with respect to
$t$, multiplying it by $r^m\rho^4u_{ttt}$, and integrating by
parts over $I$, we have \beqno
&&\int_Ir^m\rho^5u_{ttt}^2+\frac{1}{2}\frac{d}{dt}\int_I(r^m\rho^4u_{rtt}^2+mr^{m-2}\rho^4u_{tt}^2)\\&
=&2\int_Ir^m\rho^3\rho_tu_{rtt}^2-4\int_Ir^m\rho^3\rho_ru_{rtt}u_{ttt}-\int_Ir^m\rho^4u_{ttt}\Big[\rho_{tt}u_t+2\rho_tu_{tt}+\rho_{tt}uu_r+2\rho_tu_tu_r
\\&&+2\rho_tuu_{rt}+\rho u_{tt}u_r+2\rho u_tu_{rt}+\rho
uu_{rtt}+(\rho^\gamma)_{rtt}\Big]+2\int_Imr^{m-2}\rho^3\rho_tu_{tt}^2\\&&+\int_Ir^m\rho^4u_{ttt}\Big(\rho_{tt}f+2\rho_tf_t+\rho
f_{tt}\Big)\\&\le&c\int_Ir^m\rho^2u_{rtt}^2-8\int_Ir^m\rho^\frac{5}{2}u_{ttt}\rho(\sqrt{\rho})_ru_{rtt}+\frac{1}{4}\int_Ir^m\rho^5u_{ttt}^2
+c+c\int_Ir^m(f^2+f_t^2+f_{tt}^2)\\&\le&c\int_Ir^m\rho^2u_{rtt}^2+\frac{1}{2}\int_Ir^m\rho^5u_{ttt}^2
+c\int_Ir^m(f^2+f_t^2+f_{tt}^2)+c, \eeqno where we have used
(\ref{2.3}), (\ref{2.4}), (\ref{2.7}), (\ref{2.8}), Lemma
\ref{le:3.2},
 Lemma \ref{le:3.3}, Lemma \ref{le:3.5}  and Cauchy inequality.

Thus, \beqno
\int_Ir^m\rho^5u_{ttt}^2+\frac{d}{dt}\int_I(r^m\rho^4u_{rtt}^2+mr^{m-2}\rho^4u_{tt}^2)\le
c\int_Ir^m\rho^2u_{rtt}^2+c\int_Ir^m(f^2+f_t^2+f_{tt}^2)+c.\eeqno
By (\ref{1.8}), (\ref{2.2}) and $\rho_0\nabla^3{\bf g}_1\in L^2$,
we have \beqno
\int_{Q_T}r^m\rho^5u_{ttt}^2+\int_I(r^m\rho^4u_{rtt}^2+r^{m-2}\rho^4u_{tt}^2)\le
c.\eeqno The proof of  Lemma \ref{le:3.6} is complete. \qed

\begin{Lemma}\label{le:3.7}For any $0\le t\le T$, it holds
$$
\int_Ir^m(\rho^2u_{rrrt}^2+|\partial_r^5u|^2)+\int_{Q_T}r^m(\rho^3u_{rrtt}^2+\rho
u_{rrrrt}^2) \le c.
$$
\end{Lemma}
\noindent{\it Proof.} From (\ref{3.17}), (\ref{2.3}), (\ref{2.4}),
(\ref{2.7}), (\ref{2.8}), Lemmas \ref{le:3.1}-\ref{le:3.3} and
Lemma \ref{le:3.6}, we get \beqno
\int_Ir^m\rho^2u_{rrrt}^2&\le&c\int_Ir^m(f_t^2+f_r^2+f_{rt}^2)+c\\&\le&c.
\eeqno This combining (\ref{3.24}) and Lemma \ref{le:3.4} gives
$$
\int_Ir^m|\partial_r^5u|^2\le c.
$$
By (\ref{2.2}), (\ref{2.3}), (\ref{2.4}), (\ref{2.7}),
(\ref{2.8}), Lemma \ref{le:3.1}, Lemma \ref{le:3.3} and Lemma
\ref{le:3.6}, we get \beq\label{3.26}
\nonumber\int_{Q_T}r^m\rho^3u_{rrtt}^2&\le&
c\int_{Q_T}r^m(f_t^2+f_{tt}^2)+c\\&\le&c. \eeq Differentiating
(\ref{3.17}) with respect to $r$, we have
\beq\label{3.27}\nonumber
u_{rrrrt}&=&\rho_{rrt}u_t+2\rho_{rt}u_{rt}+\rho_{rr}u_{tt}+2\rho_ru_{rtt}+\rho_tu_{rrt}+\rho
u_{rrtt}+\Big[\rho_{rt}uu_r+\rho_ru_tu_r\\&&\nonumber+\rho_ruu_{rt}+\rho_tu_r^2+2\rho
u_ru_{rt}+\rho_tuu_{rr}+\rho u_tu_{rr}+\rho
uu_{rrt}+(\rho^\gamma)_{rrt}\\&&-mr^{-1}u_{rrt}+2mr^{-2}u_{rt}-\frac{2m
u_t}{r^3}\Big]_r-\rho_{rrt}f-2\rho_{rt}f_r-\rho_{rr}f_t-\rho_tf_{rr}\\&&\nonumber-2\rho_rf_{rt}-\rho
f_{rrt}. \eeq By (\ref{3.26}), (\ref{3.27}), (\ref{2.3}),
(\ref{2.4}), (\ref{2.7}), (\ref{2.8}), Lemmas
\ref{le:3.1}-\ref{le:3.3}  and Lemma \ref{le:3.5}, we obtain
\beqno \int_{Q_T}r^m\rho
u_{rrrrt}^2&\le&c\int_{Q_T}r^m\rho^2|(\sqrt{\rho})_r|^2u_{rtt}^2+c\int_{Q_T}r^m\rho^3u_{rrtt}^2\\
&&+c\int_{Q_T}r^m( f_{rr}^2+f_{rt}^2+f_{rrt}^2)+c\\&\le&c. \eeqno
The proof of  Lemma \ref{le:3.7} is complete. \qed

\begin{Lemma}\label{le:3.8}For any $0\le t\le T$, it holds
$$
\int_Ir^m(|\partial_r^5\rho|^2+|\partial_r^5(\rho^\gamma)|^2)+\int_{Q_T}r^m
|\partial_r^6u|^2 \le c.
$$
\end{Lemma}
\noindent{\it Proof.} Differentiating (\ref{3.18}) with respect to
$r$, we obtain \beq\label{3.28}
&&\nonumber\partial_r^5\rho_t+mr^{-1}\rho u_{rrrrr}+mr^{-1}\rho_r
u_{rrrr}-mr^{-2}\rho
u_{rrrr}+\Big(4mr^{-1}\rho_ru_{rrr}+6mr^{-1}\rho_{rr}u_{rr}
\\ && -\nonumber\frac{4m\rho u_{rrr}}{r^2}+\frac{12m\rho
u_{rr}}{r^3}-\frac{18m\rho
u_r}{r^4}+4mr^{-1}\rho_{rrr}u_r\Big)_r+mr^{-1}\rho_{rrrrr}u+mr^{-1}\rho_{rrrr}
u_r
\\ && -\nonumber\frac{5m\rho_{rrrr}u}{r^2}-\frac{4m\rho_{rrr}u_r}{r^2}+\frac{8m\rho_{rrr}u}{r^3}+\Big[-\frac{12m\rho_r
u_{rr}}{r^2}-\frac{12m\rho_{rr} u_r}{r^2}+\frac{24m\rho_ru_r}{r^3}
\\ && +\nonumber\frac{12m\rho_{rr}u}{r^3}-\frac{6m\rho
u_r}{r^4}-\frac{24m\rho_r u}{r^4}+\frac{24m\rho
u}{r^5}\Big]_r+6\rho_r\partial_r^5u+\rho\partial_r^6u+15\rho_{rr}u_{rrrr}
\\ && +20\rho_{rrr}u_{rrr}+15\rho_{rrrr}u_{rr}+6\partial_r^5\rho
u_r+u\partial_r^6\rho =0. \eeq Multiplying (\ref{3.28}) by
$r^m\partial_r^5\rho$, and using (\ref{2.3}), (\ref{2.4}),
(\ref{2.7}), (\ref{2.8}), Lemma \ref{le:3.1}, Lemma \ref{le:3.3},
Lemma \ref{le:3.4}  and Lemma \ref{le:3.7}, we have
\beq\label{3.29}\nonumber\frac{d}{dt}\int_Ir^m|\partial_r^5\rho|^2&\le&c\int_Ir^m|\partial_r^5\rho|^2+c\int_Ir^m|\partial_r^6u|^2
+\int_Ir^m({|\partial_r^5\rho|^2})_ru+c\\&\le&c\int_Ir^m|\partial_r^5\rho|^2+c\int_Ir^m|\partial_r^6u|^2
+c.\eeq Similarly, we get \be\label{3.30}
\frac{d}{dt}\int_Ir^m|\partial_r^5(\rho^\gamma)|^2\le
c\int_Ir^m|\partial_r^5(\rho^\gamma)|^2+c\int_Ir^m|\partial_r^6u|^2
+c. \ee Differentiating (\ref{3.22}) with respect to $r$, we get
\beq\label{3.31}\nonumber
\partial_r^6
u&=&\rho_{rrrr}u_t+4\rho_{rrr}u_{rt}+6\rho_{rr}u_{rrt}+4\rho_ru_{rrrt}+\rho
u_{rrrrt}+(\rho_ruu_r+\rho u_r^2+\rho
uu_{rr})_{rrr}\\&&\nonumber+\partial_r^5(\rho^\gamma)+\left(-mr^{-1}u_{rrrr}+4mr^{-2}u_{rrr}-12mr^{-3}u_{rr}+24mr^{-4}u_r-\frac{24m
u}{r^5}\right)_r\\&&-\rho_{rrrr}f-4\rho_{rrr}f_r-6\rho_{rr}f_{rr}-4\rho_rf_{rrr}-\rho
f_{rrrr}. \eeq (\ref{3.31}), (\ref{2.3}), (\ref{2.4}),
(\ref{2.7}), (\ref{2.8}), Lemma \ref{le:3.1}, Lemma \ref{le:3.3},
Lemma \ref{le:3.4}  and Lemma \ref{le:3.7} imply \beq
\label{3.32}\nonumber\int_Ir^m|\partial_r^6u|^2&\le&c\int_Ir^mu_{rrt}^2+c\int_Ir^mu_{rrrt}^2+c\int_Ir^m\rho
u_{rrrrt}^2+c\int_Ir^m|\partial_r^5(\rho^\gamma)|^2\\&&+c\int_Ir^m(f_{rr}^2+f_{rrr}^2+f_{rrrr}^2)+c.
\eeq It follows from (\ref{3.29}), (\ref{3.30}), (\ref{3.32}),
 (\ref{2.8}), Lemma \ref{le:3.3}, Lemma \ref{le:3.7}
and Gronwall inequality that \be\label{3.33}
\int_Ir^m[|\partial_r^5\rho|^2+|\partial_r^5(\rho^\gamma)|^2]\le
c. \ee From (\ref{3.32}), (\ref{3.33}), (\ref{2.8}), Lemma
\ref{le:3.3} and Lemma \ref{le:3.7}, we obtain
$$
\int_{Q_T}r^m|\partial_r^6u|^2\le c.
$$
The proof of   Lemma \ref{le:3.8} is complete.\qed

It follows from $(\ref{1.4})_1$, (\ref{2.6}), (\ref{2.7}),
(\ref{2.8})  and Lemmas \ref{le:3.1}-\ref{le:3.8} that
\be\label{3.34}
\|(\rho,\rho^\gamma)\|_{H_r^5}+\|(\rho_t,(\rho^\gamma)_t)\|_{H_r^4}+\|((\sqrt{\rho})_r,(\sqrt{\rho})_t)\|_{L^\infty(Q_T)}\le
c,\ \rho\ge\frac{\delta}{c}, \ee and \beq\label{3.35}
&&\nonumber\int_Ir^m(\rho^3u_{tt}^2+\rho^4u_{rtt}^2+r^{-2}\rho^4u_{tt}^2+\rho^2u_{rrrt}^2+\rho
u_{rrt}^2+u_{rt}^2+r^{-2}u_t^2+|\partial_r^5u|^2+|\partial_r^4u|^2\\&&\nonumber+u_{rrr}^2+u_{rr}^2+u_r^2+r^{-2}u^2)+\int_{Q_T}r^m(\rho
u_{tt}^2+\rho^2u_{rtt}^2+\rho^5u_{ttt}^2+\rho^3u_{rrtt}^2+\rho
u_{rrrrt}^2\\&&+u_{rrt}^2+u_{rrrt}^2+|\partial_r^6u|^2)\le c.\eeq
By (\ref{3.34}) and (\ref{3.35}), we complete the proof of Theorem
\ref{th:3.1}. \qed

\vspace{4mm}

{\bf Proof of Theorem \ref{th:1.2}}:

Since (\ref{3.34}) and (\ref{3.35}) are uniform for $b$ and
$\delta$, it suffices to prove Theorem \ref{th:1.2} for the case
$b<\infty$. We follow the strategy as the proof of Theorem
\ref{th:1.1} and use Theorem \ref{th:3.1}. After taking
$\delta\rightarrow 0$ (take subsequence if necessary), we get a
solution $(\rho,u)$ of (\ref{1.4})-(\ref{1.6}) satisfying
\beq\label{3.36}\begin{cases} (\rho,\rho^\gamma)\in
L^\infty([0,T];H^5(I)),\
((\sqrt{\rho})_r,(\sqrt{\rho})_t)\in L^\infty(Q_T),\\
(\rho_t,(\rho^\gamma)_t)\in L^\infty([0,T];H^4(I)),\  u\in
L^\infty([0,T];H^5(I))\cap L^2([0,T];H^6(I)),\\ u_t\in
L^\infty([0,T];H_0^1(I))\cap L^2([0,T];H^3(I)),\\
(\sqrt{\rho}\partial_r^2 u_t, \rho\partial_r^3
 u_t)\in L^\infty([0,T];L^2(I))\cap L^2([0,T];H^1(I)).\end{cases} \eeq It
follows from $u\in L^\infty([0,T];H^5(I))\cap L^2([0,T];H^6(I)),$
$u_t\in L^2([0,T];H^3(I)),$ (\ref{1.4})$_1$, (\ref{2.6})  and the
similar arguments as \cite{8,9} that
$$ \rho,\rho^\gamma\in C([0,T];H^5(I)).
$$

Denote $\rho({\bf x},t)=\rho(r,t),\ {\bf u}({\bf
x},t)=u(r,t)\frac{{\bf x}}{r},$ then $(\rho,{\bf u})$ is the
unique solution to (\ref{1.1})-(\ref{1.3}) with the regularities
in Theorem \ref{th:1.2}. The proof of Theorem \ref{th:1.2} is
complete. \qed

\begin{Remark}\label{re:3.1}
{\rm{ By our method, it seems that the regularities of ${\bf u}$
could not be improved to $L^\infty([0,T];D^6(\Omega))$ and
$L^2([0,T];D^8(\Omega))$, even if the initial data, ${\bf f}$ and
${\bf g_1}$ of (\ref{1.1})-(\ref{1.3}) are  smooth enough. More
precisely, based on (\ref{3.34}) and (\ref{3.35}), by using
similar arguments as in the proofs of Theorem 1.1 and
 1.2, we get the next two
 {\it a priori} estimates about $u$
 for (\ref{1.4})-(\ref{1.6}).  \beq\label{3.37}
\int_Ir^m\rho^7u_{ttt}^2+\int_{Q_T}r^m\rho^6(u_{rttt}^2+r^{-2}u_{ttt}^2)\le
c,\eeq and
\beq\label{3.38}\int_Ir^m\rho^8(u_{rttt}^2+r^{-2}u_{ttt}^2)+\int_{Q_T}r^m\rho^9u_{tttt}^2\le
c. \eeq (\ref{3.37}) and (\ref{3.38}) respectively implies
\begin{equation}\label{3.39}\int_Ir^m\rho|\partial_r^6u|^2+\int_{Q_T}r^m|\partial_r^7u|^2\le
c,\end{equation}and
\begin{equation}\label{3.40}\int_Ir^m\rho^2|\partial_r^7u|^2+\int_{Q_T}r^m\rho|\partial_r^8u|^2\le
c.\end{equation}
 But the appearance of vacuum stops us from obtaining the regularities of
${\bf u}$ in $L^\infty([0,T];D^6(\Omega))$ and
$L^2([0,T];D^8(\Omega))$ from (\ref{3.39}) and (\ref{3.40}). Is
there another way to get further regularities of the solutions,
such as $L^\infty([0,T];D^6(\Omega))-$regularity? This is an
interesting problem. }}
\end{Remark}

\vspace{6mm}

 \noindent{\bf\small Acknowledgment.} {\small

 The authors would like to thank the anonymous referees for their constructive suggestions and kindly comments. S.J. Ding was supported by the National Basic Research Program of China (973 Program) (No.2011CB808002),
the National Natural Science Foundation of China (No.11071086),
and the University Special Research Foundation for Ph.D Program
(No.20104407110002).
 L. Yao was supported by  the
National Natural Science Foundation of China    $\#$11101331 and
China Postdoctoral Science Foundation funded project
$\#$20100481359, $\#$201104676. C.J. Zhu was supported by the
National Natural Science Foundation of China $\#$10625105,
$\#$11071093, the PhD specialized grant of the Ministry of
Education of China $\#$20100144110001, and the Special Fund for
Basic Scientific Research  of Central Colleges $\#$CCNU10C01001.}

\end{document}